\documentclass[12pt]{article}

\usepackage{epsfig}
\usepackage{a4}
\usepackage{amsmath,amsthm,amssymb}
\newcommand{\be}{\begin{equation}}
\newcommand{\ee}{\end{equation}}
\newcommand{\bea}{\begin{eqnarray}}
\newcommand{\eea}{\end{eqnarray}}

\newcommand{\dofig}[2]{
  \begin{figure}[ht]
  \hbox to\linewidth{\hss
    \epsfig{bbllx=56bp,bblly=166bp,bburx=561bp,bbury=565bp,%
    file=#1.eps,width=0.6\linewidth}\hss}
  \caption[fig#1]{\label{fig#1}#2}
  \end{figure}
}

\begin{document}

\thispagestyle{empty}

\hbox to\hsize{%
  \vbox{%
%       \hbox{DRAFT}%
%       \hbox{Submitted to}%
%       \hbox{\sl ???}%
        }\hfil
  \vbox{%
%        \hbox{MPP-2004-169}%
        \hbox{May 23, 2009}%
%        \hbox{gr-qc/0412067}%
        }}

\vspace{1.5cm}
\begin{center}
{\LARGE\bf Self-similar Solutions of the \\Cubic Wave Equation\\}

\vspace{1.5 cm} \large P. Bizo\'n${}^a$, P. Breitenlohner${}^b$, D. Maison${}^b$,
A. Wasserman${}^c$\\

\vspace{1 cm}

${}^a$ {\small\sl
M. Smoluchowski Institut of Physics\\
Jagiellonian University\\
Krak\'ow, Poland\\
}

 \vspace{0.3 cm}

${}^b$ {\small\sl
Max-Planck-Institut f\"ur Physik\\
--- Werner Heisenberg Institut ---\\
Munich, Germany\\}

 \vspace{0.3 cm}

${}^c$ {\small\sl
Department of Mathematics\\
University of Michigan\\
Ann Arbor, Michigan\\}

\end{center}

\vspace{20 mm}
\begingroup \addtolength{\leftskip}{1cm} \addtolength{\rightskip}{1cm}
\begin{abstract}\noindent
We prove that the focusing cubic wave equation in three spatial dimensions has a countable family
of self-similar solutions which are smooth inside the past light cone of the singularity. These
solutions are labeled by an integer index $n$ which counts the number of oscillations of the
solution. The linearized operator around the $n$-th solution is shown to have $n+1$ negative
eigenvalues (one of which corresponds to the gauge mode) which implies that all $n>0$ solutions
are unstable. It is also shown that all $n>0$ solutions have a singularity outside the past light
cone which casts doubt on whether these solutions may participate in the Cauchy evolution, even
for non-generic initial data.
\end{abstract}

\endgroup
\newpage

\section{Introduction}
This paper is a continuation of our studies of semilinear wave equations in three spatial
dimensions with a focusing power nonlinearity
\begin{equation}\label{eqo}
\partial_{tt} u- \Delta u- u^p=0\,,\qquad p=\mbox{odd integer}\geq 3\,.
\end{equation}
In \cite{bmw} we showed that for each odd integer $p\geq 7$ equation (\ref{eqo}) has a countable
sequence of regular self-similar solutions while for $p=5$ there is no nontrivial regular
self-similar solution. This result has important consequences for the character of the threshold
of blowup for Eq.~(\ref{eqo}) \cite{bct}.

Here we consider the subcritical power $p=3$ and, as before, restrict our attention to
spherically symmetric solutions, so $u=u(t,r)$ and Eq.~(\ref{eqo}) reduces to
\begin{equation}\label{eqor}
\partial_{tt} u- \partial_{rr} u- \frac{2}{r}\partial_r u - u^3=0\,.
\end{equation}
This equation has the scaling symmetry (for each positive constant $\lambda$)
\begin{equation}\label{scal}
    u(t,r)\, \rightarrow \, u_{\lambda} (t,r) = \lambda^{-1}
    u(t/\lambda,r/\lambda)\,,
\end{equation}
so it is natural to ask if there are solutions which are invariant (modulo time translation)
under this scaling. Such solutions are called self-similar. It follows from (\ref{scal}) and the
symmetry under time translation that self-similar solutions must have the form
\begin{equation}\label{cssrho}
    u(t,r)=(T-t)^{-1} U(\rho)\,,
\end{equation}
where $T$ is a positive constant (usually referred to as the blowup time) and $\rho=r/(T-t)$ is
the so-called similarity variable ranging from zero to infinity. By definition, the self-similar
solutions are singular at the point $(T,0)$.

 Substituting the ansatz (\ref{cssrho}) into Eq.~(\ref{eqor}) we get the ordinary
 differential equation
\begin{equation}\label{eqrho}
(1-\rho^2) \frac{d^2 U}{d\rho^2} +\left(\frac{2}{\rho}-4\rho\right) \frac{dU}{d\rho}-2U+U^3=0\;.
\end{equation}
The obvious constant solution of this equation is  $U_0(\rho)=\sqrt{2}$ and the question is if
there are other nontrivial solutions $U(\rho)$ which are smooth in the interval $0\le\rho\le1$
(i.e., inside the past light cone of the point $(T,0)$). Numerical evidence for the existence of
a countable family of such solutions was given in \cite{bct} and the main goal of this paper is
to prove this fact rigorously (we note that a general variational argument for the existence of
infinitely many weak solutions of equation (\ref{eqrho}) was given before in \cite{gal}). As in
\cite{bmw}, we will present two different proofs of this result. The first proof, given in
Sects.~2 and~3, is rather explicit and exploits the conformal invariance of the cubic wave
equation is an essential way. The second proof, given in the Appendix, is more general in the
sense that it is based on 'soft' topological arguments.

 Having established the existence of self-similar solutions in Sect.~3,
 we analyze some of their
 properties in the second part of the paper. In Sect.~4 we derive some remarkable asymptotic
 scaling formulae for the solutions with many oscillations.
 Sect.~5 is devoted to the linear stability analysis.
 Finally, in Sect.~6 we show that all nonconstant solutions have a singularity
 outside the past light cone.
 We point out that for the sake of clarity of exposition the paper is written in the
 'physics' style, however its conversion to the 'epsilon-delta' style is routine and we leave it to the
 mathematically oriented reader.
\section{Dynamical system and local existence}
 In the studies of self-similar solutions
  it is convenient to use hyperbolic polar coordinates $(s,x)$ defined by
\begin{equation}\label{polar}
 T-t=e^{-s}\cosh(x)\,,\qquad r=e^{-s} \sinh(x)\,.
\end{equation}
The transformation (\ref{polar}) is a conformal transformation of the Minkowski spacetime. In the
hyperbolic coordinates the Minkowski metric reads
\begin{equation}\label{metric}
    ds^2= e^{-2s} \left(-ds^2+dx^2+\sinh^2(x)\, d\Omega^2\right)\,,
\end{equation}
where $-\infty<s<\infty, x\geq 0$ and $d\Omega^2$ is the round metric on the unit two-sphere. The
surfaces $s=const$ are hyperboloids $H^3$ with constant scalar curvature $-1$ which foliate the
interior of the past light cone of the point $(T,0)$.  Due to the conformal symmetry of
Eq.~(\ref{eqor}), the function $f(s,x)=ru(t,r)$ satisfies a simple wave equation
\begin{equation}\label{eqh}
    \partial_{ss} f - \partial_{xx} f - \frac{f^3}{\sinh^2(x)}=0\,.
\end{equation}
The self-similar solutions of Eq.~(\ref{eqor}) correspond to static solutions of Eq.~(\ref{eqh}),
i.e., solutions $f(x)$ which satisfy the ordinary differential equation (here and in the
following we denote the derivative by the prime)
\be\label{neweq} f''+\frac{f^3}{\sinh^2(x)}=0 \ee on the half-line $x\geq 0$. In particular, the
constant solution $U_0(\rho)=\sqrt{2}$ of Eq.~(\ref{eqrho}) corresponds to
$f_0(x)=\sqrt{2}\tanh(x)$. The remainder of this section and sections~3 and~4 are devoted to the
analysis of solutions of Eq.~(\ref{neweq}).

In order to obtain a dynamical system formulation we introduce
\be b(x)=f(x)-x f'(x)\;,\qquad d(x)=f'(x)\,.\ee
Then, Eq.~(\ref{neweq}) is equivalent to the system of first order equations
\begin{equation}\label{dbeq}
b' = \frac{xf^3}{\sinh^2(x)}\;, \qquad d'=-\frac{f^3}{\sinh^2(x)}\;,
\end{equation}
where $f(x)=b(x)+x d(x)$. Rewriting this system in the form
\begin{equation}\label{dbdyn}
x\left(\frac{b}{x}\right)'=
  -\frac{b}{x}+\frac{x^4}{\sinh^2(x)}\left(\frac{b}{x}+d\right)^3\;,
\qquad x d'=
  \frac{x^4}{\sinh^2(x)}\left(\frac{b}{x}+d\right)^3\;,
\end{equation}
and applying Prop.~1 of \cite{BFM}, we infer that there exists a one-parameter family
$(b/x,d)(x)$ of local solutions of Eqs.~(\ref{dbdyn}) with boundary condition $(b/x,d)(0)=(0,c)$,
analytic in $(x^2,c)$ and defined for all $c$ and $|x|<\xi(c)$ with some $\xi(c)>0$. We shall
refer to these solutions as `$c$-orbits' and denote them by $b(c,x)$ and $d(c,x)$. It follows
from the above that the $c$-orbits have the following expansion near the origin
\begin{equation}\label{exp0}
    b(c,x)=\frac{1}{3} c^3 x^3 + \mathcal{O}\left(x^5\right)\,,\qquad d(c,x)=c-\frac{1}{2} c^3
    x^2 +
    \mathcal{O}\left(x^4\right)\,.
\end{equation}

\section{Global existence}
In this section we consider the global behavior of $c$-orbits. Without loss of generality we may
assume that $c\geq 0$. First, we observe that $c$-orbits are defined for all $x\geq 0$, as
follows immediately from the existence of the
  Lyapunov function
\begin{equation}\label{ljap}
G=2d^2+\frac{f^4}{\sinh^2(x)}\;, \qquad G'=-2\coth(x)\frac{f^4}{\sinh^2(x)}\le0\;.
\end{equation}
Second, $c$-orbits have simple asymptotic behavior since by (\ref{ljap}) $f/\sqrt{\sinh(x)}$ is
bounded when $x\to\infty$, hence the right hand sides of Eqs.~(\ref{dbeq}) are integrable, and
therefore $b(x)$ and $d(x)$ have finite limits $B$ and $D$ for $x\to\infty$. Moreover, the rapid
decrease of the right hand sides of Eqs.~(\ref{dbeq}) implies that the limits $B(c)$ and $D(c)$
are continuous functions of $c$.

  It follows from the above that $f(c,x)\sim B(c)+ D(c) x$ for $x\rightarrow \infty$, hence solutions of
  Eq.~(\ref{neweq}) corresponding to $c$-orbits
  are regular at infinity   if and only if $D(c)=0$.
  Now, we will show that there
  is an infinite sequence $c_n$ with $n=0$, 1, $2,\ldots$ such that $D(c_n)=0$. The
  corresponding globally regular solutions are
  characterized by $b_n=B(c_n)$ and behave asymptotically as
  \begin{equation}\label{basym}
    f(x)=b_n - b_n^3 e^{-2x} + \mathcal{O}(e^{-4x})\,.
  \end{equation}
 Numerically one finds a sequence of such
solutions (see Fig.~\ref{figgb}) with  $n=0,1,2,\ldots$ zeros and parameters $b_n^2\sim
c_n\sim(n+1)^3$ (see Table~1 and Sect.~4).
Note that for the globally regular solutions the conserved energy functional associated with
Eq.~(\ref{eqh}),
\begin{equation}\label{energy}
    E(f)=\frac{1}{2} \int_0^{\infty} \left(f_s^2+f_x^2-\frac{f^4}{2\sinh^2(x)}\right) \, dx\,,
\end{equation}
is finite and, as indicated by numerics, monotonically increasing with $n$ (see Table~1).
 \dofig{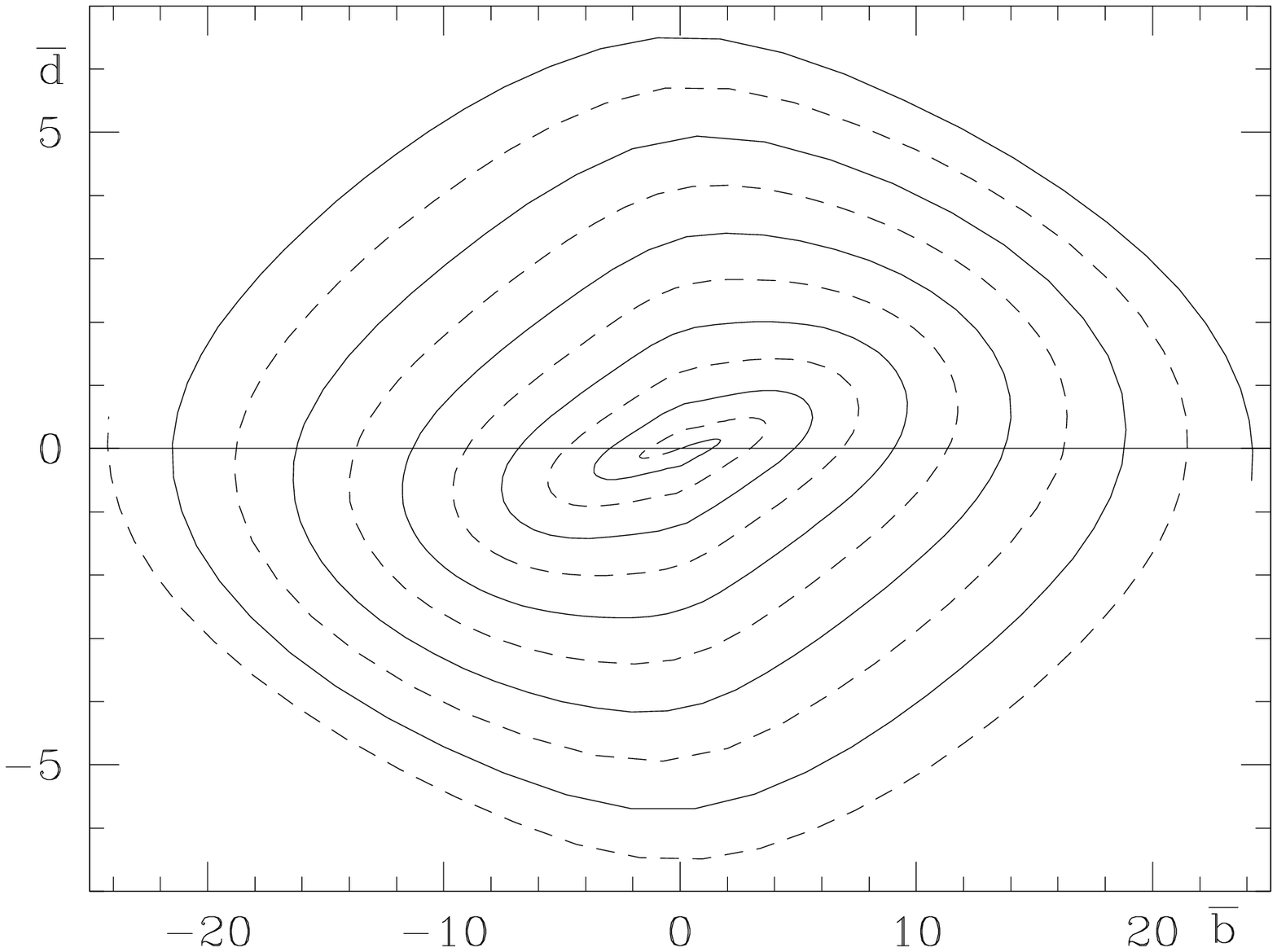}{$\bar d(c)=D/\sigma$ vs.\ $\bar
b(c)=(B+8D)/\sigma$ where $\sigma=((B+8D)^2+D^2)^{1/6}$ for positive (solid) and negative
(dashed) values of $c$.}
\begin{table}[h]
\begin{center}
\begin{tabular*}
{0.82\textwidth}{@{\extracolsep{\fill}}|c|c|c|c|c|c|}
  \hline
  $n$ & $c_n$ & $b_n$ & $E(f_n)$ & $c_n$ (theory) & $b_n$ (theory)\\
  \hline
  0 & $\sqrt{2}$ & $\sqrt{2}$ & 1/3 & 1.630626 & 1.467029 \\
  1 & $$9.616283 & -3.578348 & 4.62810 & 9.991135 & -3.631358\\
  2 & 30.13927 & 6.315947 &  21.5429 & 30.681145 & 6.363520\\
  3 & 68.58242 & -9.519976 & 64.8053& 69.292246  & -9.563216\\
  4 & 130.5379 & 13.13018 & 153.071 & 131.41603  & 13.170001\\
  5 & 221.5967 & -17.10516 & 309.116 & 222.64408  &-17.142226\\
  6 & 347.3277 & 21.41418 & 556.682 & 348.56798 & 21.448919 \\
  \hline
\end{tabular*}
  \caption{The parameters of the first few solutions $f_n$ generated numerically
  and their comparison with the asymptotic formulae (\ref{cn}) and (\ref{bn}) .}
  \end{center}
\end{table}

Let us introduce the phase function
\be \phi(x)=\arctan\left(\frac{b(x)}{d(x)}\right)\;, \qquad {\rm with} \quad
\phi'=\frac{f^4}{\sinh^2(x)(b^2+d^2)}\ge0\;. \ee
For $c$-orbits we normalize $\phi(c,x)$ by the condition $\phi(c,0)=0$. Then
$\phi(c,x)=(i-1/2)\pi$ at the $i^{\rm th}$ extremum of $f(c,x)$. Furthermore, the limit
$\Phi(c)=\phi(c,\infty)$ is a continuous function of $c$ and $\Phi(c_n)=(n+1/2)\pi$ for a regular
solution $f_n(x)$ with $n$ zeros.

Integrating Eqs.~(\ref{dbeq}) for $c\approx0$ yields
\be \phi(c,x)=c^2\int_0^x\frac{\xi^4d\xi}{\sinh^2(\xi)}+\mathcal{O}(c^4)\;, \qquad
\Phi(c)=c^2\frac{\pi^4}{30}+\mathcal{O}(c^4)\;. \ee
To find the behavior of the phase function for large $c$ we rescale the variables
\begin{equation}\label{rescale}
    F(y)=f(x)\,,\qquad y=c x\,,
\end{equation}
so that Eq.~(\ref{neweq}) becomes \be\label{neweq2} F''+\frac{F^3}{c^2 \sinh^2(y/c)}=0\,. \ee For
$c\rightarrow\infty$ we get the limiting equation
\be\label{neweq3} F''+\frac{F^3}{y^2}=0\,, \qquad F(y)=y-\frac{1}{6}y^3+\mathcal{O}(y^5)\,, \ee
whose solution is oscillatory. In terms of the original variables this implies that
  $\lim_{c\to\infty}\phi(c,x)=\infty$ for any finite $x>0$, and hence
$\lim_{c\to\infty}\Phi(c)=\infty$. Therefore, for each $n\ge0$ there exists at least
one\footnote{Numerics indicates that there for each $n$ there is exactly one $c_n$.} value $c_n$
such that $\Phi(c_n)=(n+1/2)\pi$. This concludes the proof of existence of a countable family of
regular self-similar solutions of the cubic wave equation.
\section{Asymptotic formula for $b_n$ and $c_n$}
In this section, we use the technique of matched asymptotic expansions to derive the asymptotic
scaling formulae for the parameters of solutions $f_n(x)$ in the large $n$ limit.
 The solutions
with many zeros are approximately periodic with a modulated amplitude. In order to extract the
periodic part we factorize $f$ in the form
\be f(x)=a(x)v(t(x))\;, \ee
with suitable functions $a(x)$ and $t(x)$. Plugging this ansatz into Eq.~(\ref{neweq}) we get
(denoting $t$-derivatives by a dot)
\be\label{ansatz} at'^2\ddot v+(at''+2a't')\dot v+a''v+\frac{a^3v^3}{\sinh^2(x)}=0\;, \ee
and impose the conditions
\bea\label{cond}
&&\sinh(x)t'=a\;,\\
&&2a't'+at''=0\;. \eea
Differentiating the first equation we get
\be \frac{a'}{a}=\frac{t''}{t'}+\coth(x)\;, \ee
which we use to eliminate $a$ from the second one. Thus we obtain
\be \frac{t''}{t'}=-\frac{2}{3}\coth(x)\;, \ee
and by integration (suppressing an irrelevant integration constant)
\be t'=\sinh^{-2/3}(x)\;, \ee
and finally
\be t(x)=\int_0^x \frac{d\xi}{\sinh^{2/3}(\xi)}\;. \ee
The length $T$ of the $t$ interval corresponding to $0\leq x<\infty$ is
\be T=\int_0^{\infty}
\frac{d\xi}{\sinh^{2/3}(\xi)}=\frac{1}{2}B\left(\frac{1}{6},\frac{1}{3}\right)
\approx4.20654632\;. \ee
From Eq.~(\ref{cond}) we get
\be a(x)=\sinh^{1/3}(x)\;. \ee
Using the expressions for $a(x)$ and $t(x)$ in Eq.~(\ref{ansatz}) yields
\be\label{cublin} \ddot v+v^3+hv=0\;, \ee
with
\be h(x)=\frac{3\sinh^2(x)-2\cosh^2(x)}{9\sinh^{2/3}(x)}\;. \ee
The behavior of $t(x)$ for $x\to 0$ is $t(x)\to 3x^{1/3}$, implying $h\to -2/t^2$ for $x\to 0$.
Introducing $\bar t=T-t$ we find $\bar t(x)\to \frac{3}{2}\cosh^{-2/3}(x)$ for $x\to\infty$ and
consequently $h\to 1/4\bar t^2$.  The linear term $hv$ can be neglected except near $t=0$ resp.\
$\bar t=0$, where it dominates the cubic term. From the boundary conditions for $f(x)$ one
obtains $v(t)\to\frac{c}{9}t^2$ for $t\to 0$ resp. $v(\bar t)\to \sqrt{\frac{2}{3}}b\bar t^{1/2}$
for $\bar t\to 0$.

In order to extract the leading behaviour for large $b$ resp.\ $c$, we rescale $t\to t/c^{1/3}$
resp.\ $\bar t\to \bar t/b^{2/3}$. Furthermore we rescale $v$ by $v\to c^{1/3}v$ resp.\ $v\to
b^{2/3}v$. In the limit $b,c\to\infty$, neglecting nonleading terms in $h$ we obtain the
equations
\be\label{asy0} \ddot v+v^3-\frac{2}{t^2}v=0\;, \ee
resp.
\be\label{asy1} \ddot v+v^3+\frac{1}{4\bar t^2}v=0\;, \ee
for the rescaled variables. The rescaled boundary conditions are $v(t)\to \frac{1}{9}t^2$ for
$t\to 0$ resp.\ $v(\bar t)\to \sqrt{\frac{2}{3}}\bar t^{1/2}$ for $\bar t\to 0$. Numerically one
finds that the solutions of Eqs.~(\ref{asy0},\ref{asy1}) with these boundary conditions converge
very quickly to solutions of
\be\label{cubic} \ddot v+v^3=0\;, \ee
with amplitudes $A_0$ resp.\ $A_1$. Numerically one finds $A_0\approx0.90247851$ and
$A_1\approx0.82273965$. After rescaling $v(t)\to v(t/A)/A$ both solutions tend to the solution
$F_1(t)$ of Eq.~(\ref{cubic}) with the normalized amplitude and the corresponding period
\be \tau=4\sqrt2\int_0^1\frac{dz}{(1-z^4)^{1/2}}=
     \sqrt{2}B(\frac{1}{4},\frac{1}{2})\approx7.41629871\;.
\ee
We fix the phase of $F_1$ such that $F_1(0)=0$. Then the rescaled $v$'s tend to $F_1$ with some
phaseshifts $\theta_0$ resp. $\theta_1$, i.e.
\be v(t/A)/A\to F_1(t+\theta)\; \qquad {\rm for}\quad t\to\infty\;. \ee
Numerically one finds $\theta_0\approx-1.6225533$ and $\theta_1\approx0.8623512$.  After shifting
the two solutions with their respective $\theta$'s they coincide asymptotically.

Supposing we have a regular solution $f_n(x)$ with $n\gg1$ zeros, the parameters $b_n$ and $c_n$
must have been chosen so that the corresponding solutions of Eq.~(\ref{cublin}) starting at
either end of the interval $0\leq x<\infty$ match at some intermediate point. Using the discussed
asymptotics of the rescaled solutions we obtain two conditions matching amplitudes and phases.
Matching the amplitudes we get
\be\label{amp} c_n^{1/3}A_0=b_n^{2/3}A_1\;. \ee
To match the phases we assume that we match a solution starting at $x=0$ with $m$ zeros with one
from $x=\infty$ with $n-m+1$ at their last zero. The corresponding $t$ intervals $t_0$ and $\bar
t_1$ must add up to the total $t$ interval $T$. Taking into account the rescalings and the
phaseshifts we obtain the condition
\be \frac{c_n^{-1/3}}{A_0}\Bigl(m\frac{\tau}{2}-\theta_0\Bigr)
+\frac{b_n^{-2/3}}{A_1}\Bigl((n-m+1)\frac{\tau}{2}-\theta_1\Bigr)=T\;. \ee
Making use of Eq.~(\ref{amp}) this can be rewritten as
\be \frac{(n+1)\frac{\tau}{2}-(\theta_0+\theta_1)}{c_n^{1/3}A_0}=T\;, \ee
which is independent of $m$ and hence from the matching point as required for consistency. Thus
we obtain the desired asymptotic formula for $n\gg1$
\be c_n=\Biggl(\frac{(n+1)\frac{\tau}{2}-(\theta_0+\theta_1)}{A_0T}\Biggr)^3\;. \ee
Putting in numbers we get
\be\label{cn} c_n\approx\Biggl(\frac{3.70814935\,(n+1)+0.7602022}{3.7963177}\Biggr)^3\;, \ee
together with
\be\label{bn} b_n^2=\Biggl(\frac{A_0}{A_1}\Biggr)^3 c_n
  \approx1.3198462\,c_n\;.
\ee
\section{Linear stability analysis}
The role of self-similar solutions in the Cauchy evolution depends crucially on their stability
under small perturbations. To analyze this issue, in this section we determine the spectrum of
the linearized operator around the self-similar solutions $f_n(x)$. Substituting the ansatz
$f(s,x)=f_n(x)+w(s,x)$ into Eq. (\ref{eqh}) and linearizing, we obtain the linear wave equation
with a potential
\begin{equation}\label{eqw}
    \partial_{ss} w - \partial_{xx} w + V_n(x) w=0\,,\qquad V_n(x)=-\frac{3 f_n^2}{\sinh^2(x)}\,.
\end{equation}
Since $f_n(x)\sim c_n x$ for $x\rightarrow 0$ and $f_n(\infty)=b_n$, the potential $V_n(x)$ is
everywhere bounded and decays exponentially at infinity.

Separating time, $w(s,x)=e^{i k s} \xi(x)$, we get the eigenvalue problem
\begin{equation}\label{eqhlin}
 L_n \xi = k^2 \xi\,, \qquad L_n=-\frac{d^2}{dx^2} + V_n(x)\,.
\end{equation}
For each $n$, the operator $L_n$ is self-adjoint on $\mathcal{D}(L_n)=\{\xi\in L_2[0,\infty),
\xi(0)=0\}$ and has a continuous spectrum $k^2\geq 0$. The discrete spectrum depends on $n$. More
precisely, we claim that $L_n$ has exactly $n+1$ negative eigenvalues. To see this, note that
Eq.~(\ref{eqh}) is invariant under the following transformation
\begin{equation}\label{transhyp}
s \rightarrow s+\frac{1}{2} \ln\left(1+2\alpha \cosh(x) e^s+\alpha^2 e^{2s}\right)\,,\quad x
\rightarrow \tanh^{-1}\left(\frac{\sinh(x)}{\cosh(x)+\alpha e^s}\right)\,,
\end{equation}
which is nothing else but the time translation $t \rightarrow t+\alpha$, expressed in hyperbolic
coordinates. Hence, each time-independent solution $f(x)$ gives rise to the one-parameter family
of time-dependent solutions
\begin{equation}\label{orbith}
    f_{\alpha}(s,x)=f\left(\tanh^{-1}\left(\frac{\sinh(x)}{\cosh(x)+\alpha e^s}\right)\right)\,.
\end{equation}
The  perturbation $\delta f$ generated by this symmetry (which we shall refer to as the gauge
mode) is tangent to the orbit (\ref{orbith}) at $\alpha=0$, that is
\begin{equation}\label{gaugemode}
    \delta f= \frac{\partial f_{\alpha}(s,x)}{\partial \alpha}\vert_{\alpha=0} = \sinh(x) f'(x) \,e^s\,.
\end{equation}
Thus, for each $n$ the operator $L_n$  has the eigenvalue $k^2=-1$ associated with the
eigenfunction $\xi^{(n)}(x)=\sinh(x) f_n'(x)$. Since by construction the solution $f_n(x)$ has
$n$ oscillations, it follows that the eigenfunction $\xi^{(n)}(x)$ has $n$ nodes, which in turn
implies by the Sturm oscillation theorem that there are exactly $n$ eigenvalues below $k^2=-1$.
This means that the solution $f_n(x)$ has at least $n$ unstable modes (the gauge mode does not
count as a genuine instability). Numerics indicate that there are no eigenvalues in the interval
$-1<k^2<0$, so we claim that the above phrase 'at least $n$' can be replaced by 'exactly $n$',
however we can prove this claim only for the perturbations of the ground state solution
$f_0(x)=\sqrt{2}\tanh(x)$. In this case
\begin{equation}\label{pt}
    V_0(x)=-\frac{6}{\cosh^2(x)}
\end{equation}
is the P\"oschl-Teller potential for which the whole spectrum is known explicitly. In particular,
the gauge mode $\xi^{(0)}(x)=\sinh(x)/\cosh^2(x)$ is the only eigenfunction.
\section{Behavior of solutions outside the light cone}
The hyperbolic coordinates (\ref{polar}) cover only the interior of the past light cone, hence in
order to see what happens outside the light cone we need to return to the similarity variable
$\rho$ and  Eq.~(\ref{eqrho}). The results of Sect.~3 imply that Eq.~(\ref{eqrho}) has infinitely
many solutions which are smooth on the closed interval $0\leq \rho \leq 1$ and behave as
\begin{equation}\label{local0}
U(\rho) = c + \frac{c}{6}(2 - c^2) \rho^2 + \mathcal{O}(\rho^4) \quad \text{for} \quad \rho
\rightarrow 0\,,
\end{equation}
and
\begin{equation}\label{local1}
U(\rho) = b + \frac{b}{2} (b^2 -  2) (\rho -1) +\mathcal{O}((\rho-1)^2) \quad \text{for} \quad
\rho \rightarrow 1\,.
\end{equation}
We will refer to solutions satisfying the boundary condition (\ref{local1}) as to  '$b$-orbits'.
Without loss of generality we assume that $b\geq 0$. Now, we will show that there are no
$b$-orbits which are smooth for all $\rho\geq 0$. The proof consists of two steps. In the first
step we show that  all $b$-orbits with $0<b<\sqrt{2}$ become singular at $\rho=0$. This will
imply that the smooth solutions constructed in Section~3 must have $b_n>\sqrt{2}$ for all $n\geq
1$. In the second step we show that all $b$-orbits with $b>\sqrt{2}$ become singular at some
$\rho>1$.
 \vskip 0.2cm \noindent \emph{Step~1.}  Consider a $b$-orbit with
$0<b<\sqrt{2}$ and assume that it exists for all $\rho\leq 1$ and is smooth at $\rho=0$. Let us
define the function
\begin{equation}\label{h}
h(\rho)=-\frac{U'(\rho)}{U(\rho)}\,.
\end{equation}
Using equation (\ref{eqrho}) we get
\begin{equation}\label{hp}
h'(\rho)=\frac{(\rho-\rho^3)U'^2+(2-4\rho^2)UU'-\rho U^2 (2-U^2)}{\rho(1-\rho^2) U^2}\,.
\end{equation}
 It follows from (\ref{local1}) that $h(1)=\dfrac{1}{2}(2-b^2)>0$ and $h'(1)=\dfrac{1}{8}(2-b^2)(b^2-4)<0$.
 We will show that $h'(\rho)<0$ for all $\rho>0$ and therefore $\lim_{\rho\rightarrow 0^+}h(\rho)\geq
 h(1)>0$. But smoothness at $\rho=0$ requires that $h(0)=0$. This contradiction will prove Step~1.
 To show that $h'(\rho)<0$ we show equivalently that the numerator of the fraction on the right hand side
  of Eq.(\ref{hp})
 \begin{equation}
 n(\rho)=(\rho-\rho^3)U'^2+(2-4\rho^2)UU'-\rho U^2(2-U^2)
\end{equation}
is non-positive.
 From (\ref{local1}) we have $n(1)=0$ and
 $n'(1)=\dfrac{1}{4}b^2(4-b^2)(2-b^2)>0$, hence $n(\rho)$ is negative in the left neighborhood of $\rho=1$.
  If $n(\rho)$ is negative for all $\rho<1$ then we are
 done. Thus, let us suppose that there is a $\rho=R$ such that
 $n(R)=0$ and $n(R)<0$ for $R<\rho<1$. We will show that this is impossible
 because $n'(R)>0$. Note that equation
 \begin{equation}\label{n(R)}
n(R)=(R-R^3)U'^2+(2-4R^2)UU'-R U^2
 (2-U^2)=0\,,
\end{equation}
viewed as a formal quadratic equation for $U'$, can be satisfied only if the discriminant
(divided by a positive factor $4U^2$ for convenience)
\begin{equation}
 \Delta=1-R^2(1-R^2)(2+U^2)
 \end{equation}
 is non-negative.
 A calculation yields
 \begin{equation}\label{gp}
n'(R)=\frac{2U}{R(1-R^2)} \left[-(\Delta+R^2) U'+R^3 U(U^2-2)\right]\,,
 \end{equation}
where terms involving $U'^2$ were eliminated using Eq.(\ref{n(R)}). If $U(R)\geq \sqrt{2}$ then
the right hand side of Eq.(\ref{gp}) is manifestly positive. The case $0<U(R)<\sqrt{2}$ requires
more work. In this case Eq.(\ref{n(R)})
  has a single negative root
 \begin{equation}
 U'(R)=\frac{\left(1-2R^2+\sqrt{\Delta}\right) U}
 {R(R^2-1)}\,.
 \end{equation}
Substituting this value into (\ref{gp}) we get $n'(R)=(positive factor)\cdot N(U,R)$, where
\begin{equation}\label{G}
N(U,R)=(\Delta+R^2)(1-2R^2+\sqrt{\Delta})+(1-R^2)R^4(U^2-2)\,.
\end{equation}
One can check that in the rectangle $0\leq U \leq \sqrt{2}$, $0\leq R\leq 1$ the function
$N(U,R)$ has no critical points so its minima and maxima occur on the boundary. It is easy to
verify that $N(U,R)\geq 0$ on all sides of the rectangle, thus $N(U,R)>0$ inside the rectangle.
This implies that $n'(R)>0$ and completes the proof of Step~1.

\vskip 0.2cm \noindent \emph{Step~2:} Consider a $b$-orbit with $b>\sqrt{2}$ and assume that it
exists for all $\rho\geq 1$. It follows immediately from Eq.~(\ref{eqrho}) that $U(\rho)$ is
monotone increasing.
 Next, let
\begin{equation}\label{function}
    g(\rho)=\rho^4 U(\rho)U'(\rho)-\frac{\rho^3}{6} (U^2(\rho)-2)\,.
\end{equation}
From (\ref{local1}) we get $g(1)=b(b^2-2)/3>0$. We claim that $g(\rho)>0$ for all $\rho\geq 1$.
To show this let us compute $g'(\rho)$ at a point where $g(\rho)=0$. After a straightforward
calculation we get
\begin{equation}
    g'(g=0)=\frac{\rho^2(U^2-2)}{18(\rho^2-1)}\left((17U-9)\rho^2+U+3\right)\,.
\end{equation}
Since the right hand side of this equation is manifestly positive for $\rho>1$, $g(\rho)$ cannot
cross zero from above, hence $g(\rho)>0$ for all $\rho\geq 1$, as claimed. Thus
\begin{equation}\label{heq}
    \frac{U'}{U^2-2} > \frac{1}{6\rho}\,,
\end{equation}
which after integration from $1$ to $\rho$ yields a contradiction
\begin{equation}\label{hint}
    \frac{U(\rho)-\sqrt{2}}{U(\rho)+\sqrt{2}} \geq \frac{b-\sqrt{2}}{b+\sqrt{2}}
    \, \exp\left(\frac{\ln(\rho)}{3\sqrt{2}}\right)\,,
\end{equation}
which completes the proof of Step~2.
 \vskip 0.2cm
 It follows from the above analysis that each self-similar solution $U_n(\rho)$ with $n\geq 1$ is
singular at some $\rho_n>1$. This fact suggests that these solutions do not participate in the
Cauchy evolution of smooth initial data which in turn corroborates the conjecture that the
constant solution $U_0=\sqrt{2}$ is a universal attractor for blowup solutions (see \cite{mz} for
what is known rigorously and \cite{bz} for recent numerical evidence) .

\subsubsection*{Acknowledgements} PB is grateful to Frank Merle for helpful remarks.
The work of PB was supported in part by the MNII grants NN202 079235 and 189/6.PRUE/2007/7.

\section*{Appendix}

\noindent We present here an alternative  'soft' topological proof of existence of infinitely
many smooth self-similar solutions of Eq.(5). \vskip 0.2cm We already know that given any $c$
there is a unique smooth solution $U(\rho,c)$ of Eq.(5) satisfying $U(0,c)=c$ defined for all
$0\leq \rho <1$. Similarly, there is a unique smooth solution $U(\rho,b)$ satisfying $U(1,b)=b$
defined for all $0< \rho \leq 1$.
 \vskip 0.2cm
We express these solutions in terms of polar coordinates, that is we define
\begin{equation}\label{rtheta}
    r(\rho,c)=\sqrt{U(\rho,c)^2+U'(\rho,c)^2}\,, \quad \theta(\rho,c)=\arctan\left(\frac{
    U'(\rho,c)}{U(\rho,c)}\right)\,,
\end{equation}
and similarly
\begin{equation}\label{Rbeta}
    R(\rho,b)=\sqrt{U(\rho,b)^2+U'(\rho,b)^2}\,, \quad \beta(\rho,b)=\arctan\left(\frac{
    U'(\rho,b)}{U(\rho,b)}\right)\,.
\end{equation}
Let $\rho_0=\sqrt{2/3}$. Since the region $\{(\rho,c)| 0<\rho\leq \rho_0,c> 0\}$ (resp.
$\{(\rho,b)| \rho_0\leq \rho\leq 1,b> 0\}$) is simply connected, the angle $\theta(\rho,c)$
(resp. $\beta(\rho,b )$) is defined unambiguously once we specify its value at any point in the
domain. We set $\theta(0,1)=0$, hence $\theta(0,c)=0$ for all $c>0$. Similarly, we set
$\beta(1,\sqrt{2})=0$; then $-\pi/2<\beta(1,b)<\pi/2$ for all $b>0$.
Next, we define maps
\begin{equation}
    \Phi: R_{+}\ni c \rightarrow \Phi(c)=(\theta(\rho_0,c),R(\rho_0,c))\in R^2_{+}\,,
\end{equation}
and
\begin{equation}
    \Psi_k: R_{+}\ni
   b \rightarrow \Psi_k(b)=(\beta(\rho_0,b)-2k\pi,R(\rho_0,b))\in R^2_{+}\,.
\end{equation}
Note that if $\Psi_k(b)=\Phi(c)$ for some $b$ and $c$, then we have a solution defined on the
whole interval $0\leq \rho\leq 1$ in the nodal class with index $n=2k$ (according to our
terminology from Section~3).

 \vskip 0.2cm \noindent Lemma 1.
 $\lim_{c\rightarrow 0}
r(\rho_0,c)=0$ and  $\lim_{b\rightarrow 0} R(\rho_0,b)=0$.\\
\emph{Proof:} Follows immediately from continuous dependence on initial conditions.
\vskip 0.2cm \noindent Lemma 2. $\lim_{c\rightarrow \infty}
\theta(\rho_0,c)=-\infty$ and  $\lim_{b\rightarrow \infty} \beta(\rho_0,b)=\infty$.\\
\emph{Proof:} Follows from the asymptotic analysis given in Section~4.
 \vskip 0.2cm \noindent Lemma 3. For any positive $b$ and $c$ we have
 $\theta(\rho,c)<\pi/2$ and $\beta(\rho,b)>-\pi/2$.\\
 \emph{Proof:} We have $\theta(0,c)=0$ and if $\theta(\rho,c)=\pi/2$, theni
 $U(\rho,c)=0,U'(\rho,c)>0$ so $\theta'(\rho,c)<0$, contradiction. Similarly, we have
 $\beta(1,b)>-\pi/2$ and if $\beta(\rho,b)=-\pi/2$ then $U(\rho,b)=0,U'(\rho,b)<0$ so
 $\beta'(\rho,b)>0$, contradiction.
  \vskip 0.2cm \noindent Lemma 4. If $0<c<2$, then $-\pi/2<\theta(\rho_0,c)<\pi/2$ and
  similarly, if $0<b<2$, then $-\pi/2<\beta(\rho_0,b)<\pi/2$.\\
  \emph{Proof:} We define the function
  \begin{equation}\label{funH}
  H(\rho)=\frac{1}{2}(1-\rho^2)U'^2-U^2+\frac{1}{4}U^4\,.
  \end{equation}
We have $H'(\rho)=(3\rho-2/\rho)U'^2$ so $H(\rho)$ decreases on $(0,\rho_0]$ and increases on
$[\rho_0,1]$.  If $0<c<2$ and $\rho\leq \rho_0$ then $H(\rho,c)<H(0,c)<0$, hence $U(\rho,c)>0$
(since $H\geq 0$ if $U=0$). Similarly, if $0<b<2$ and $\rho\geq \rho_0$ then
$H(\rho,b)<H(1,b)<0$, hence $U(\rho,b)>0$.
 \vskip 0.2cm \noindent
Now we are ready to prove: \vskip 0.2cm \noindent \textbf{Theorem.}
 For any positive integer $n$ there exist parameters $(c_n,b_n)$ such that the corresponding
 solution
 $U(\rho,c_n)=U(\rho,b_n)$ is in the $n^{th}$ nodal class. \vskip 0.2cm \noindent \emph{Proof:}
If $n=2k$ then by Lemmas 2 and 3, for any integer $k\geq 1$ we may choose $c_R>c_L>2$ (resp.
$b_R>b_L>2$) such that $\theta(\rho_0,c_L)=-\pi/2$ and $\theta(\rho_0,c_R)=-(2k+1)\pi$ (resp.
$\beta(\rho_0,b_L)=\pi/2$ and $\beta(\rho_0,b_R)=(2k+1)\pi$). Let $c_R$ (resp. $b_R$) be the
smallest such $c$ (resp. $b$). Then $-\pi/2>\theta(\rho_0,c)>-(2k+1)\pi$ for $b_L<b<b_R$.Then
$\pi/2<\beta(\rho_0,b)<(2k+1)\pi$ for $b_L<b<b_R$. Next, we choose $m$ and $M$ such that
$m<r(\rho_0,c)<M$ for $c_L<c<c_R$ and $m<R(\rho_0,b)<M$ for $b_L<b<b_R$. Finally, by Lemma~1 we
choose $\tilde c<c_L$ (resp. $\tilde b<b_L$) such that $r(\rho_0,\tilde c))=m$ (resp.
$R(\rho_0,\tilde b)=m$) and let $\tilde c$ (resp. $\tilde b$) be the largest such $c$ (resp.
$b$). Let $\Omega$ be the rectangle with vertices
$(-(2k+1)\pi,m)$,$(-(2k+1)\pi,M)$,$(\pi,m)$,$(\pi,M)$. The ordered points $A=\Psi(\tilde b),
B=\Phi(\tilde c), C=\Psi(b_R), D=\Phi(c_R)$ lie on the boundary of $\Omega$, thus it follows from
elementary topology that the curve $\{\Phi(c)|\tilde c\leq c\leq c_R\}$ from $B$ to $D$  and the
curve $\{\Psi_k(b)|\tilde b\leq b\leq b_R\}$ from $A$ to $C$ must intersect.

If $n= 2k + 1$ we can repeat the above argument with $b < 0$ making the obvious modifications in
Lemma~3 ($\beta(\rho,b) > \pi/2$) and Lemma~4 ($\pi/2 < \beta(\rho,b) < 3\pi/2$).

\end{document}